\documentclass[leqno,a4paper]{article}
\usepackage{amsfonts}
\usepackage{graphics,graphicx,amssymb}
\usepackage{amsmath,amsthm}

\numberwithin{equation}{section} \numberwithin{figure}{section}
\newtheorem{theorem}{Theorem}[section]
\newtheorem{corollary}[theorem]{Corollary}

\newfont{\got}{eufm9 scaled 1095}

\newfont{\bbb}{msbm10}

\newcommand{\norm}[1]{\left\Vert#1\right\Vert}

\begin{document}

\title{Lie groups as four-dimensional special complex manifolds with Norden metric}
\author{Marta Teofilova}
\date{}
\maketitle
\begin{abstract}
An example of a four-dimensional special complex manifold with
Norden metric of constant holomorphic sectional curvature is
constructed via a two-parametric family of solvable Lie algebras.
The curvature properties of the obtained manifold are studied.
Necessary and sufficient conditions for the manifold to be isotropic
K\"ahlerian are given.

\noindent\textbf{2000 Mathematics Subject Classification:} 53C15,
53C50.

\noindent\textbf{Keywords:} almost complex manifold, Norden metric,
Lie group, Lie algebra.
\end{abstract}

\section{Preliminaries}
Let $(M,J,g)$ be a $2n$-dimensional almost complex manifold with
Norden metric, i.e. $J$ is an almost complex structure and $g$ is a
metric on $M$ such that:
\begin{equation}
J^{2}x=-x,\qquad g(Jx,Jy)=-g(x,y),\qquad x,y\in \mbox{\got X}(M).
\label{1.1}
\end{equation}
The associated metric $\widetilde{g}$ of $g$ on $M$, given by
$\widetilde{g} (x,y)=g(x,Jy)$, is a Norden metric, too. Both metrics
are necessarily neutral, i.e. of signature $(n,n)$.

If $\nabla $ is the Levi-Civita connection of $g$, the tensor field
$F$ of type $(0,3)$ is defined by
\begin{equation}
F(x,y,z)=g\left((\nabla _{x}J)y,z\right) \label{F}
\end{equation}
and has the following symmetries
\begin{equation}\label{Fs}
F(x,y,z)=F(x,z,y)=F(x,Jy,Jz).
\end{equation}

Let $\left\{ e_{i}\right\} $ ($i=1,2,\ldots ,2n$) be an arbitrary
basis of $ T_{p}M$ at a point $p$ of $M$. The components of the
inverse matrix of $g$ are denoted by $g^{ij}$ with respect to the
basis $\left\{ e_{i}\right\} $. The Lie 1-forms $\theta $ and
$\theta^{\ast}$ associated with $F$ are defined by, respectively
\begin{equation}\label{1-3}
\theta (x)=g^{ij}F(e_{i},e_{j},x), \qquad \theta^{\ast}=\theta \circ
J.
\end{equation}

The Nijenhuis tensor field $N$ for $J$ is given by
\begin{equation}\label{N}
N(x,y)=[Jx,Jy]-[x,y]-J[Jx,y]-J[x,Jy].
\end{equation}
It is known \cite{N-N} that the almost complex structure is complex
iff it is integrable, i.e. $N=0$.

A classification of the almost complex manifolds with Norden metric
is introduced in \cite{Ga-Bo}, where eight classes of these
manifolds are characterized according to the properties of $F$. The
three basic classes: $\mathcal{W}_1$, $\mathcal{W}_2$ of \emph{the
special complex manifolds with Norden metric} and $\mathcal{W}_3$ of
\emph{the quasi-K\"ahler manifolds with Norden metric} are given as
follows:
\begin{equation}
\begin{array}{l}
\mathcal{W}_{1}:F(x,y,z)=\frac{1}{2n}\left[ g(x,y)\theta
(z)+g(x,z)\theta (y)\right.
\smallskip \\
\qquad \qquad \qquad \qquad \quad \left. +g(x,Jy)\theta
(Jz)+g(x,Jz)\theta (Jy)\right] ;\medskip \\
\mathcal{W}_{2}:F(x,y,Jz)+F(y,z,Jx)+F(z,x,Jy)=0,\quad \theta =0 \hspace{0.07in}\Leftrightarrow\hspace{0.07in} N=0,\quad \theta =0;\medskip \\
\mathcal{W}_{3}:F(x,y,z)+F(y,z,x)+F(z,x,y)=0.
\end{array}\label{1.5}
\end{equation}
The class $\mathcal{W}_0$ of \emph{the K\"ahler manifolds with
Norden metric} is defined by $F=0$ and is contained in each of the
other classes.

Let $R$ be the curvature tensor of $\nabla $, i.e. $R(x,y)z=\nabla
_{x}\nabla _{y}z-\nabla _{y}\nabla _{x}z-\nabla _{\left[ x,y \right]
}z$. The corresponding (0,4)-type tensor is defined by
$R(x,y,z,u)=g\left( R(x,y)z,u\right)$. The Ricci tensor $\rho$ and
the scalar curvatures $\tau$ and $ \tau^{\ast}$ are given by:
\begin{equation}
\begin{array}{c}
\rho(y,z)=g^{ij}R(e_{i},y,z,e_{j}),\qquad
\tau=g^{ij}\rho(e_{i},e_{j}),\qquad \tau^{\ast}=g^{ij}\rho
(e_{i},Je_{j}).
\end{array}
\label{tau}
\end{equation}

A tensor of type (0,4) is said to be \emph{curvature-like} if it has
the properties of $R$. Let $S$ be a symmetric (0,2)-tensor. We
consider the following curvature-like tensors:
\begin{equation}\label{psi}
\begin{array}{l}
\psi_{1}(S)(x,y,z,u) =g(y,z)S(x,u)-g(x,z)S(y,u) \smallskip\\
\phantom{\psi_{1}(S)(x,y,z,u)}+  g(x,u)S(y,z) -
g(y,u)S(x,z),\smallskip\\
\pi_{1}=\frac{1}{2}\psi_{1}(g), \quad \pi_2(x,y,z,u)= g(y,Jz)g(x,Ju)
- g(x,Jz)g(y,Ju).
\end{array}
\end{equation}

It is known that on a pseudo-Riemannian manifold $M$ ($\dim M=2n
\geq 4$) the conformal invariant Weyl tensor has the form
\begin{equation}\label{Weyl}
W(R)=R-\frac{1}{2(n-1)}\big
\{\psi_{1}(\rho)-\frac{\tau}{2n-1}\pi_{1}\big \}.
\end{equation}

Let $\alpha =\left\{ x,y\right\} $ be a non-degenerate $2$-plane
spanned by the vectors $x,y\in T_{p}M$, $p\in M$. The sectional
curvature of $\alpha $ is given by
\begin{equation}
k (\alpha ;p)=\frac{R(x,y,y,x)}{\pi _{1}(x,y,y,x)}. \label{sec}
\end{equation}
We consider the following basic sectional curvatures in $T_{p}M$
with respect to the structures $J$ and $g$: \emph{holomorphic
sectional curvatures} if $J\alpha =\alpha $ and \emph{totally real
sectional curvatures} if $J\alpha \perp \alpha $ with respect to
$g$.

The square norm of $\nabla J$ is defined by $\norm{\nabla J}
^{2}=g^{ij}g^{kl}g\left( (\nabla _{e_{i}}J)e_{k},(\nabla
_{e_{j}}J)e_{l}\right)$. Then, by (\ref{F}) we get
\begin{equation}\label{n1}
\norm{\nabla J} ^{2} = g^{ij}g^{kl}g^{pq}F_{ikp}F_{jlq},
\end{equation}
where $F_{ikp}=F(e_i,e_k,e_p)$.

An almost complex manifold with Norden metric satisfying the
condition $\norm{\nabla J}^{2}=0$ is called an \emph{isotropic
K\"{a}hler manifold with Norden metric} \cite{Mek-Man}.

\section{Almost complex manifolds with Norden metric of constant holomorphic sectional curvature}

In this section we obtain a relation between the vanishing of the
holomorphic sectional curvature and the vanishing of $\norm{\nabla
J}^2$ on $\mathcal{W}_2$-manifolds and $\mathcal{W}_3$-manifolds
with Norden metric.

In \cite{Dje-Gri} it is proved the following

\noindent\textbf{Theorem A. (\cite{Dje-Gri})} \emph{An almost
complex manifold with Norden metric is of pointwise constant
holomorphic sectional curvature if and only if}
\begin{equation}\label{R-h}
\begin{array}{l}
3\{R(x,y,z,u)+R(x,y,Jz,Ju)+R(Jx,Jy,z,u)+R(Jx,Jy,Jz,Ju)\}\smallskip\\
-R(Jy,Jz,x,u)+R(Jx,Jz,y,u)-R(y,z,Jx,Ju)+R(x,z,Jy,Ju)\smallskip\\
-R(Jx,z,y,Ju)+R(Jy,z,x,Ju)-R(x,Jz,Jy,u)+R(y,Jz,Jx,u)\smallskip\\
= 8 H\{\pi_1+\pi_2\}
\end{array}
\end{equation}
\emph{for some} $H\in FM$ \emph{and all} $x,y,z,u\in
\mathfrak{X}(M)$. \emph{In this case} $H(p)$ \emph{is the
holomorphic sectional curvature of all holomorphic non-degenerate
2-planes in} $T_pM$, $p\in M$.{\ \rule{0.4em}{0.4em}}

\smallskip

Taking into account (\ref{tau}) and (\ref{psi}), the total trace of
(\ref{R-h}) implies
\begin{equation}\label{H}
H(p) = \frac{1}{4n^2}(\tau + \tau^{\ast\ast}),
\end{equation}
where $\tau^{\ast\ast}=g^{il}g^{jk}R(e_i,e_j,Je_k,Je_l)$.

In \cite{Teo4} we have proved that on a $\mathcal{W}_2$-manifold it
is valid
\begin{equation}\label{11}
\norm{\nabla J}^2 = 2(\tau + \tau^{\ast\ast}),
\end{equation}
and in \cite{Mek-Man} it is proved that on a
$\mathcal{W}_3$-manifold
\begin{equation}\label{12}
\norm{\nabla J}^2 = -2(\tau + \tau^{\ast\ast}).
\end{equation}
Then, by Theorem A, (\ref{H}), (\ref{11}) and (\ref{12}) we obtain
\begin{theorem}\label{thh}
Let $(M,J,g)$ be an almost complex manifold with Norden metric of
pointwise constant holomorphic sectional curvature $H(p)$, $p\in M$.
Then
\begin{description}
\item[(i)] $\norm{\nabla J}^2 = 8n^2 H(p)$ \ if
$(M,J,g)\in \mathcal{W}_2$\emph{;}
\item[(ii)] $\norm{\nabla J}^2 =- 8n^2 H(p)$ \ if
$(M,J,g)\in \mathcal{W}_3$.{\ \rule{0.4em}{0.4em}}
\end{description}
\end{theorem}
Theorem \ref{thh} implies
\begin{corollary}\label{th1}
Let $(M,J,g)$ be a $\mathcal{W}_2$-manifold or
$\mathcal{W}_3$-manifold of pointwise constant holomorphic sectional
curvature $H(p)$, $p\in M$. Then, $(M,J,g)$ is isotropic K\"ahlerian
iff $H(p)=0$.
\end{corollary}

In the next section we construct an example of a
$\mathcal{W}_2$-manifold of constant holomorphic sectional
curvature.

\section{Lie groups as four-dimensional $\mathcal{W}_2$-manifolds}

Let ${\mbox{\got g}}$ be a real 4-dimensional Lie algebra
corresponding to a real connected Lie group $G$. If $\left\{
X_{1},X_{2},X_{3},X_{4}\right\}$ is a basis of left invariant vector
fields on $G$ and $[X_{i},X_{j}]=C_{ij}^{k}X_{k}$ ($i,j,k=1,2,3,4$)
then the structural constants $C_{ij}^k$ satisfy the
anti-commutativity condition $C_{ij}^k=-C_{ji}^k$ and the Jacobi
identity
$C_{ij}^{k}C_{ks}^{l}+C_{js}^{k}C_{ki}^{l}+C_{si}^{k}C_{kj}^{l}=0$.

We define an almost complex structure $J$ and a compatible metric
$g$ on $G$ by the conditions, respectively:
\begin{equation}
JX_{1}=X_{3},\quad JX_{2}=X_{4},\quad JX_{3}=-X_{1},\quad
JX_{4}=-X_{2}, \label{J-cond}
\end{equation}
\begin{equation}
\begin{array}{l}
g(X_{1},X_{1})=g(X_{2},X_{2})=-g(X_{3},X_{3})=-g(X_{4},X_{4})=1,\medskip \\
g(X_{i},X_{j})=0,\quad i\neq j,\quad i,j=1,2,3,4.
\end{array}
\label{metric}
\end{equation}
Because of (\ref{1.1}), (\ref{J-cond}) and (\ref{metric}) $g$ is a
Norden metric. Thus, $(G,J,g)$ is a 4-dimensional almost complex
manifold with Norden metric.

From (\ref{metric}) it follows that the well-known Levi-Civita
identity for $g$ takes the form
\begin{equation}\label{LC}
2g(\nabla_{X_i}X_j, X_k) = g([X_i,X_j],X_k) + g([X_k,X_i],X_j) +
g([X_k,X_j],X_i).
\end{equation}
Let us denote $F_{ijk}=F(X_i,X_j,X_k)$. Then, by (\ref{F}) and
(\ref{LC}) we have
\begin{equation}\label{Fijk}
\begin{array}{l}
2F_{ijk}= g\big([X_{i},JX_{j}]-J[X_{i},X_{j}],X_{k}\big) +
g\big(J[X_{k},X_{i}]-[JX_{k},X_{i}],X_{j}\big) \medskip\\
\qquad\qquad\qquad\qquad + g\big([X_{k},JX_{j}]-[JX_{k},X_{j}],X_{i}
\big).
\end{array}
\end{equation}

According to (\ref{1.5}) to construct an example of a
$\mathcal{W}_2$-manifold we need to find sufficient conditions for
the Nijenhuis tensor $N$ and the Lie 1-form $\theta$ to vanish on
$\mathfrak{g}$.

By (\ref{F}), (\ref{N}), (\ref{metric}) and (\ref{Fijk}) we compute
the essential components $N_{ij}^k$ ($N(X_i,X_j)=N_{ij}^kX_k$) of
$N$ and $\theta_i=\theta(X_i)$ of $\theta$, respectively, as
follows:
\begin{equation}\label{N1}
\begin{array}{l}
N_{12}^1 = C_{34}^1 - C_{12}^1 - C_{23}^3 + C_{14}^3,\qquad \theta_1
= 2C_{13}^1 - C_{12}^4 + C_{14}^2 + C_{23}^2 - C_{34}^4,\medskip\\
N_{12}^2=C_{34}^2 - C_{12}^2 - C_{23}^4 + C_{14}^4,\qquad \theta_2 =
2C_{24}^2 + C_{12}^3 + C_{14}^1+ C_{23}^1 + C_{34}^3,\medskip\\
N_{12}^3 = C_{34}^3 - C_{12}^3 + C_{23}^1 - C_{14}^1,\qquad \theta_3
= 2C_{13}^3 + C_{12}^2 + C_{14}^4 + C_{23}^4 + C_{34}^2,\medskip\\
N_{12}^4=C_{34}^4 - C_{12}^4 + C_{23}^2 - C_{14}^2,\qquad \theta_4 =
2C_{24}^4 - C_{12}^1+ C_{14}^3 + C_{23}^3 - C_{34}^1.
\end{array}
\end{equation}
Then, (\ref{1.5}) and (\ref{N1}) imply
\begin{theorem}
Let $(G,J,g)$ be a 4-dimensional almost complex manifold with Norden
metric defined by (\ref{J-cond}) and (\ref{metric}). Then, $(G,J,g)$
is a $\mathcal{W}_2$-manifold iff for the Lie algebra $\mathfrak{g}$
of $G$ are valid the conditions\emph{:}
\begin{equation}\label{ex1}
\begin{array}{l}
C_{13}^1 = C_{12}^4 - C_{23}^2 = C_{34}^4 - C_{14}^2,\qquad
C_{13}^3=-\left(C_{12}^2 +
C_{23}^4\right)=-\left(C_{14}^4+C_{34}^2\right),\medskip\\
C_{24}^4 = C_{12}^1 - C_{14}^3 = C_{34}^1 - C_{23}^3, \qquad
C_{24}^2 =
-\left(C_{12}^3+C_{14}^1\right)=-\left(C_{23}^1+C_{34}^3\right),
\end{array}
\end{equation}
where $C_{ij}^k$ ($i,j,k=1,2,3,4$) satisfy the Jacodi identity.{\
\rule{0.4em}{0.4em}}
\end{theorem}
One solution to (\ref{ex1}) and the Jacobi identity is the
2-parametric family of solvable Lie algebras $\mathfrak{g}$ given by
\begin{equation}\label{g1}
\mathfrak{g}:
\begin{array}{l}
[X_1,X_2]=\lambda X_1 - \lambda X_2,\qquad [X_2,X_3]=\mu X_1 +
\lambda X_4, \medskip\\
\lbrack X_1,X_3]=\mu X_2 + \lambda X_4, \qquad [X_2,X_4]=\mu X_1 +
\lambda
X_3, \medskip\\
\lbrack X_1,X_4]=\mu X_2 + \lambda X_3, \qquad [X_3,X_4]=- \mu X_3 +
\mu X_4,\qquad \lambda,\mu \in \mathbb{R}.
\end{array}
\end{equation}

Let us study the curvature properties of the
$\mathcal{W}_2$-manifold $(G,J,g)$, where the Lie algebra
$\mathfrak{g}$ of $G$ is defined by (\ref{g1}).

By (\ref{metric}), (\ref{LC}) and (\ref{g1}) we obtain the
components of the Levi-Civita connection:
\begin{equation}\label{LC1}
\begin{array}{ll}
\nabla_{X_1}X_2 = \lambda X_1 + \mu(X_3 + X_4),& \quad
\nabla_{X_2}X_1
= \lambda X_2 + \mu(X_3+X_4),\medskip\\
\nabla_{X_3}X_4 = -\lambda(X_1+X_2)-\mu X_3,& \quad
\nabla_{X_4}X_3=-\lambda(X_1+X_2)-\mu X_4,\medskip\\
\nabla_{X_1}X_1 = -\lambda X_2,\quad \nabla_{X_2}X_2=-\lambda X_1,&
\quad \nabla_{X_3}X_3=\mu X_4,\quad \nabla_{X_4}X_4=\mu
X_3,\medskip\\
\nabla_{X_1}X_3=\nabla_{X_1}X_4=\mu X_2,&\quad
\nabla_{X_2}X_3=\nabla_{X_2}X_4=\mu X_1,\medskip\\
\nabla_{X_3}X_1=\nabla_{X_3}X_2=-\lambda X_4,& \quad
\nabla_{X_4}X_1=\nabla_{X_4}X_2=-\lambda X_3.
\end{array}
\end{equation}
Taking into account (\ref{Fijk}) and (\ref{g1}) we compute the
essential non-zero components of $F$:
\begin{equation}\label{F1}
\begin{array}{l}
F_{114}=-F_{214}=F_{312}=\frac{1}{2}F_{322}=\frac{1}{2}F_{411}=F_{412}=-\lambda,\medskip\\
F_{112}=\frac{1}{2}F_{122}=\frac{1}{2}F_{211}=F_{212}=-F_{314}=F_{414}=\mu.
\end{array}
\end{equation}
The other non-zero components of $F$ are obtained from (\ref{Fs}).

By (\ref{n1}) and (\ref{F1}) for the square norm of $\nabla J$ we
get
\begin{equation}\label{n2}
\norm{\nabla J}^2 = -32(\lambda^2 - \mu^2).
\end{equation}

Further, we obtain the essential non-zero components
$R_{ijks}=R(X_i,X_j,X_k,X_s)$ of the curvature tensor $R$ as
follows:
\begin{equation}\label{R1}
\begin{array}{l}
-\frac{1}{2}R_{1221}=-R_{1341}=-R_{2342}=R_{3123}=\frac{1}{2}R_{3443}=R_{4124}=\lambda^2
+ \mu^2,\medskip\\
R_{1331}=R_{1441}=R_{2332}=R_{2442}=-R_{1324}=-R_{1423}=\lambda^2 - \mu^2,\medskip\\
R_{1231}=R_{1241}=R_{2132}=R_{2142}\smallskip\\=-R_{3143}=-R_{3243}=-R_{4134}=-R_{4234}=2\lambda\mu.
\end{array}
\end{equation}
Then, by (\ref{tau}) and (\ref{R1}) we get the components
$\rho_{ij}=\rho(X_i,X_j)$ of the Ricci tensor and the values of the
scalar curvatures $\tau$ and $\tau^{\ast}$:
\begin{equation}\label{tau1}
\begin{array}{l}
\rho_{11}=\rho_{22}=-4\lambda^2,\qquad\qquad\quad
\rho_{33}=\rho_{44}=-4\mu^2,\medskip\\
\rho_{12}=\rho_{34}=-2(\lambda^2+\mu^2),\qquad
\rho_{13}=\rho_{14}=\rho_{23}=\rho_{24}=4\lambda\mu,\medskip\\
\tau= -8(\lambda^2 - \mu^2),\qquad\qquad\quad\hspace{0.05in}
\tau^{\ast}=16\lambda\mu.
\end{array}
\end{equation}

Let us consider the characteristic 2-planes $\alpha_{ij}$ spanned by
the basic vectors $\{X_i,X_j\}$: totally real 2-planes -
$\alpha_{12}$, $\alpha_{14}$, $\alpha_{23}$, $\alpha_{34}$ and
holomorphic 2-planes - $\alpha_{13}$, $\alpha_{24}$. By (\ref{sec})
and (\ref{R1}) for the sectional curvatures of the holomorphic
2-planes we obtain
\begin{equation}\label{s1}
k(\alpha_{13})=k(\alpha_{24})=-(\lambda^2 - \mu^2).
\end{equation}
Then it is valid
\begin{theorem}
The manifold $(G,J,g)$ is of constant holomorphic sectional
curvature.
\end{theorem}

Using (\ref{Weyl}), (\ref{R1}) and (\ref{tau1}) for the essential
non-zero components\\ $W_{ijks}=W(X_i,X_j,X_k,X_s)$ of the Weyl
tensor $W$ we get:
\begin{equation}\label{W1}
\begin{array}{l}
\frac{1}{2}W_{1221}=W_{1331}=W_{1441}=W_{2332}=W_{2442}=\frac{1}{2}W_{3443}\medskip\\
=-\frac{1}{3}W_{1324}=-\frac{1}{3}W_{1423}=\frac{1}{3}(\lambda^2-\mu^2).
\end{array}
\end{equation}

Finally, by (\ref{Weyl}), (\ref{n2}), (\ref{tau1}), (\ref{s1}) and
(\ref{W1}) we establish the truthfulness of
\begin{theorem}
The following conditions are equivalent\emph{:}
\begin{description}
\item[(i)] $(G,J,g)$ is isotropic K\"{a}hlerian\emph{;}

\item[(ii)] $|\lambda|=|\mu|$\emph{;}

\item[(iii)] $\tau=0$\emph{;}

\item[(iv)]  $(G,J,g)$ is of zero holomorphic sectional curvature\emph{;}

\item[(v)] the Weyl tensor vanishes.

\item[(vi)] $R=\frac{1}{2}\psi _{1}(\rho
)$.

\end{description}
\end{theorem}

\begin{tabbing}
  Faculty of Mathematics and Informatics,\\
  University of Plovdiv,\\
  236 Bulgaria Blvd.,
  Plovdiv 4003, Bulgaria.\\
  e-mail:\ \texttt{marta@uni-plovdiv.bg}
  \end{tabbing}


\begin{thebibliography}{9}


\bibitem{Dje-Gri} G.~Djelepov, K.~Gribachev, \emph{Generalized $B$-manifolds of constant
holomorphic sectional curvature}, Plovdiv Univ. Sci. Works -- Math.
\textbf{23}(1) (1985), 125--131.


\bibitem{Ga-Bo} G.~Ganchev, A.~Borisov, \emph{Note on the almost complex
manifolds with a Norden metric}, Compt. Rend. Acad. Bulg. Sci.
\textbf{39}(5) (1986), 31--34.



\bibitem{Mek-Man} D.~Mekerov, M.~Manev, \emph{On the geometry of Quasi-K\"{a}hler
manifolds with Norden metric}, Nihonkai Math. J. \textbf{16}(2)
(2005), 89--93.


\bibitem{N-N} A.~Newlander, L.~Niremberg, \emph{Complex analytic coordinates
in almost complex manifolds}, Ann.~Math. \textbf{65} (1957),
391--404.

\bibitem{Teo4} M.~Teofilova, \emph{Lie groups as four-dimensional conformal K\"{a}hler manifolds
with Norden metric}, In:Topics of Contemporary Differential
Geometry, Complex Analysis and Mathematical Physics, eds. S. Dimiev
and K. Sekigawa, World Sci. Publ., Hackensack, NJ (2007), 319--326.

\end{thebibliography}
\end{document}